\documentclass[reqno, a4paper,12pt]{amsart}

\usepackage[capposition=top]{floatrow}
\usepackage{euscript,eufrak,verbatim}
\usepackage{graphicx}
\usepackage[usenames]{color}
\usepackage[colorlinks,linkcolor=red,anchorcolor=blue,citecolor=blue]{hyperref}
\usepackage{amsmath, mathtools}
\usepackage{mathtools}
\usepackage{amsthm}
\usepackage[all]{xy}
\usepackage{amssymb} 
\usepackage{bm}

\usepackage[all]{xy}

\usepackage{mathrsfs}
\usepackage{amscd}

\usepackage{enumitem}
\usepackage[numbers]{natbib}

\usepackage{mathptmx}
\usepackage[T1]{fontenc}

\usepackage{lipsum}

\makeatletter
\numberwithin{equation}{section}

\theoremstyle{cuplain}
\newtheorem{main theorem}{Main Theorem}
\newtheorem{theorem}{Theorem}[section]

\newtheorem{conjecture}[theorem]{Conjecture}
\newtheorem{corollary}[theorem]{Corollary}
\newtheorem{proposition}[theorem]{Proposition}

\theoremstyle{definition}
\newtheorem{definition}[theorem]{Definition}

\newtheorem{example}[theorem]{Example}
\newtheorem{notation}[theorem]{Notation}

\newtheoremstyle{break}
  {\topsep}{\topsep}%
  {\itshape}{}%
  {\bfseries}{}%
  {\newline}{}%
\theoremstyle{break}

\newtheoremstyle{break}
  {\topsep}{\topsep}%
  {\normalshape}{}%
  {\bfseries}{}%
  {\newline}{}%

\newtheorem{breakremark}[theorem]{Remark}

\numberwithin{equation}{section}
\newcommand{\spa}{\hspace{1pt}}
\newcommand{\vep}{\varepsilon}
\newcommand{\flo}{\mathscr}
\newcommand{\diam}{\mathrm{diam}}

\DeclarePairedDelimiter{\abs}{\lvert}{\rvert}
\DeclarePairedDelimiter{\ceil}{\lceil}{\rceil}
\newcommand{\norm}[1]{\left\lVert#1\right\rVert}
\newcommand{\setcond}{\hspace{2pt} \middle| \hspace{2pt}}

\begin{document}


\newcommand\titlelowercase[1]{\texorpdfstring{\lowercase{#1}}{#1}}

\font\mathptmx=cmr12 at 12pt


\title[\fontsize{13}{12}\mathptmx {\it{I\titlelowercase{mproved dimension theory of sofic self-affine fractals}}}]{\Huge I\titlelowercase{mproved dimension theory of} \protect{\\[7pt]} \titlelowercase{sofic self-affine fractals}}


\author[\fontsize{13}{12}\mathptmx {\it{N\titlelowercase{ima} A\titlelowercase{libabaei}}}]{\fontsize{13}{12}\mathptmx Nima Alibabaei}



\maketitle

\begin{abstract}
We establish a combinatorial expression for the Hausdorff dimension of a given self-affine fractal in any Euclidean space. This formula includes the extension of the work by Kenyon and Peres (1996) on planar sofic sets and yields an exact value for the dimension of certain sofic sets in $\mathbb{R}^3$ or higher. We also calculate the Minkowski dimension of sofic sets and establish a sufficient and presumably necessary condition for planar sofic sets to have the same Minkowski and Hausdorff dimension. The condition can be regarded as a generalization of the classical result for Bedford-McMullen carpets.
\end{abstract}

\section{Introduction} \label{section: introduction}

\subsection{Overview}

\begin{figure}[h!]
\includegraphics[width=15.8cm]{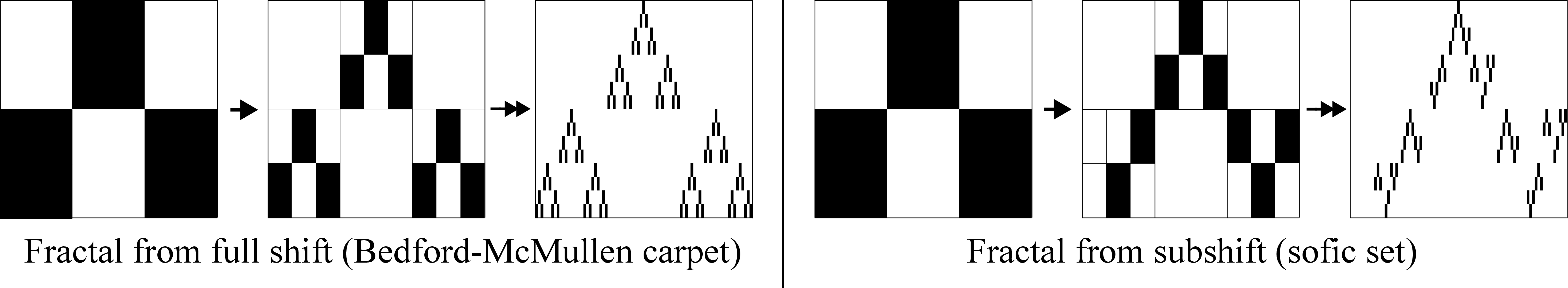}
\vspace{-10pt}%
\caption{The first few generations of self-affine fractals} \label{figure: self-affine fractals}
\end{figure}

The self-affine fractals similar to Bedford-McMullen carpets have been studied thoroughly. Figure \ref{figure: self-affine fractals} exhibits the construction of these sets. Fix integers $1 < m_1 < m_2$, and partition the square into $m_1 m_2$ numbers of congruent rectangles ($m_1 = 2$ and $m_2 = 3$ in Figure \ref{figure: self-affine fractals}). We pick some of them and continue this procedure indefinitely, where we choose the rectangles to keep at each step according to the initial choice. The generalization of these sets to arbitrary Euclidean dimensions is known as {\it{\textbf{self-affine sponges}}}, and the Hausdorff dimension and also the Minkowski dimension are all known to us \cite{Kenyon--Peres}. Behind these sets are the symbolic dynamical systems that are full shifts, an infinite Cartesian product of finite symbols: $\{0, 1, \ldots, n\}^{\mathbb{N}}$. The (left) shift-invariant property of these symbolic dynamical systems yields the affine-invariant property (or ``fractal-ness'') of self-affine sponges.

Let us consider a subset of a full shift that is shift-invariant (subshift). From it, we can define a fractal in Euclidean space similarly to how the self-affine sponges are defined. One interesting and large class of subshifts is the {\it{\textbf{sofic systems}}}, which are defined from finite digraphs. Their corresponding fractals are called the {\it{\textbf{sofic sets}}} (Figure \ref{figure: self-affine fractals}). Sofic sets in $\mathbb{R}^2$ were studied in \cite{Kenyon--Peres: sofic}, and some ingenious methods were introduced to calculate the Hausdorff dimension for various examples. Despite the intriguing demonstrations in \cite{Kenyon--Peres: sofic}, no further attempts were made to investigate the sofic sets in $\mathbb{R}^3$ or higher. The reason is that the method in \cite{Kenyon--Peres: sofic} was inherently limited to the $2$-dimensional case.

This paper aims to investigate the dimensions of such self-affine fractals. Our key result is the combinatorial formula for the Hausdorff dimension of self-affine fractals defined from arbitrary subshifts (Theorem \ref{theorem: combinatorial expression of Hausdorff dimension}). Its proof is motivated by the recent result in \cite{Alibabaei} concerning weighted topological entropy. With this result, the problem of calculating the Hausdorff dimension of self-affine fractals is reduced to counting the word complexity of the corresponding symbolic dynamical system. This result is proven in Euclidean space of arbitrary dimension. Thus, it contains the generalization of the work by Kenyon and Peres.

However, it is difficult to give an actual value using the formula in most cases. We will give an exact estimate for some sofic sets in Section \ref{section: dimension theory of sofic sets}.

Furthermore, we establish a sufficient (and presumably necessary) condition for sofic sets in $\mathbb{R}^2$ to have the same Hausdorff and Minkowski dimension (Corollary \ref{corollary: uniform complexity implies no dimension hiatus}). It is well-known that for a Bedford-McMullen carpet to have the same Hausdorff and Minkowski dimension, the number of rectangles in each row in the initial pattern must be equal for all rows. In the case of sofic sets, we see that if the number of each symbol in each row in the initial pattern is equal for all rows, then the two dimensions coincide. All the previous examples considered by Olivier \cite{Olivier} and Kenyon-Peres \cite{Kenyon--Peres: sofic} with no dimension hiatus satisfy the said condition. We prove that a more general condition is, in fact, sufficient for no dimension hiatus and conjecture that it is necessary.

\subsection{Results}

In this section, we will take a brief look at our results. The construction of sofic sets is as follows. Let $D$ be a set of labels; for example, $D = \{0, 1\} \times \{0, 1, 2\}$ in Figure \ref{figure: self-affine fractals}. We cut the square into congruent rectangular pieces and label them using elements in $D$. Each rectangle is then divided into smaller rectangles, which are labeled in the same manner. We continue this procedure inductively, then an element in $D^{\mathbb{N}}$ designates a point in the square, although the correspondence is not injective. Now consider a digraph $G$ with $D$-labeled edges. Then, an infinite path in $G$ is an element in $D^{\mathbb{N}}$; let $S \subset D^{\mathbb{N}}$ be the collection of these paths, called a sofic system. By the correspondence above, $S$ determines a self-affine fractal known as sofic set in Euclidean space as in Figure \ref{figure: self-affine fractals}. (For a rigorous definition, see Section \ref{section: results}.)

\begin{figure}[h!]
\includegraphics[width=15cm]{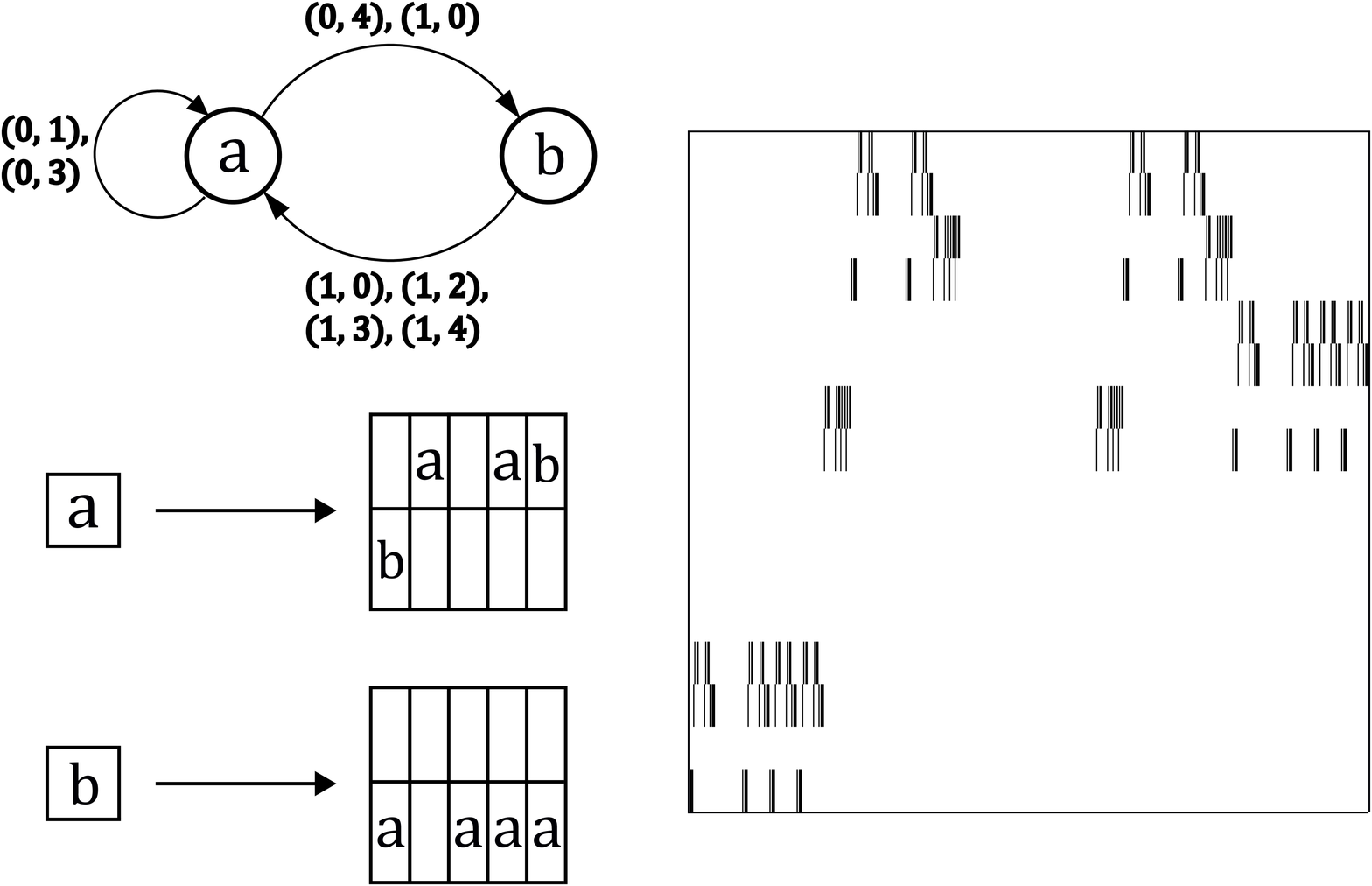}
\caption{A sofic set in $\mathbb{R}^2$ with identical Minkowski and Hausdorff dimension} \label{figure: coincidence of definition}
\end{figure}

Our main result is the combinatorial expression for the Hausdorff dimension of arbitrary self-affine fractals, including sofic sets (Theorem \ref{theorem: combinatorial expression of Hausdorff dimension}). In the case of sofic sets, we also give the formula for the Minkowski dimension (Proposition \ref{proposition: Minkowski dimension}). From these results, we establish a sufficient (and presumably necessary) condition for sofic sets in $\mathbb{R}^2$ to have the same Hausdorff and Minkowski dimension (Corollary \ref{corollary: uniform complexity implies no dimension hiatus}). For example, the sofic set in Figure \ref{figure: coincidence of definition} is constructed by starting from a square with label $a$ and following the substitution rule depicted there. Surprisingly, it has the same Minkowski and Hausdorff dimension; the dimension is $1 + \log_5{2}$.

\begin{figure}[h!]
\includegraphics[width=10cm]{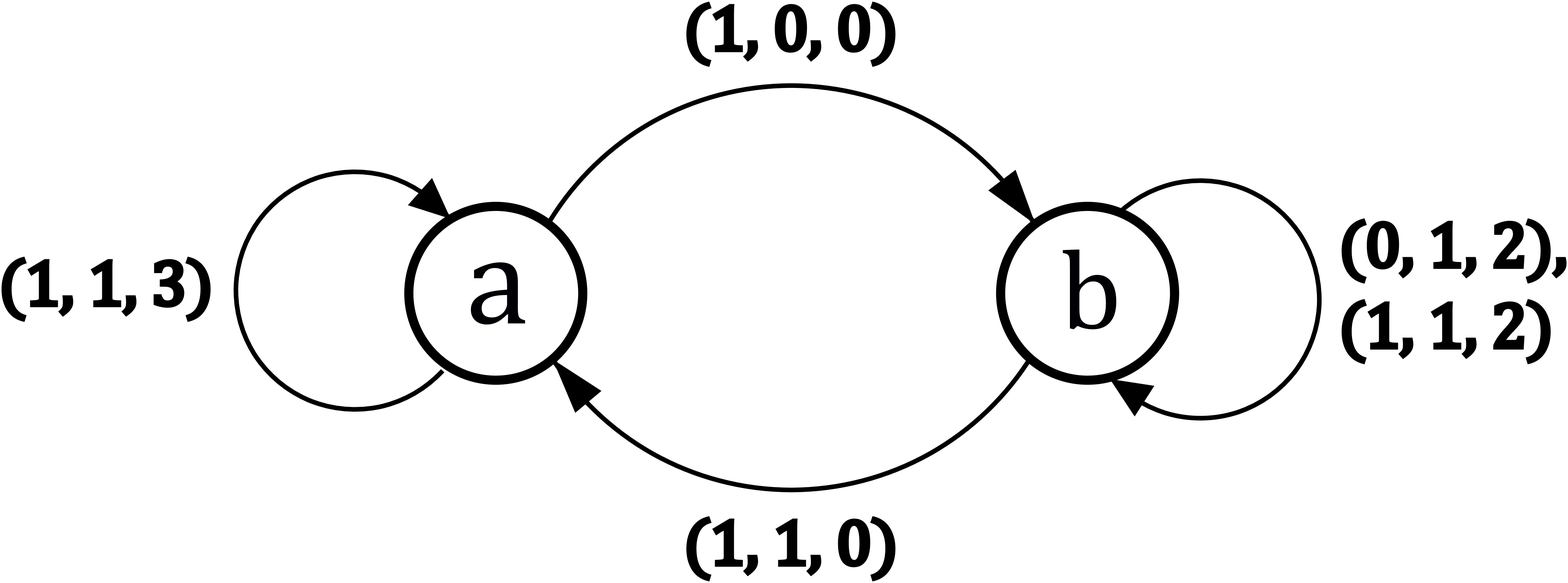}
\caption{Digraph for a sofic set in $\mathbb{R}^3$ with a calculation involving infinite matrix} \label{figure: digraph for sofic set with calculation using infinite matrix}
\end{figure}

Furthermore, we will show that we can calculate the Hausdorff dimension of some sofic sets in $\mathbb{R}^3$. For example, consider the sofic set $X$ in $\mathbb{R}^3$ defined from the digraph in Figure \ref{figure: digraph for sofic set with calculation using infinite matrix} with $D = \{0, 1\} \times \{0, 1, 2\} \times \{0, 1, 2, 3\}$. For a natural number $q$, let
\begin{equation*}
b_q =
\norm{
{\begin{pmatrix}
0 & 1 & 2^{\log_4{3}} & 3^{\log_4{3}} & 4^{\log_4{3}} & \cdots \\
1 & 0 & 0 & \cdots \\
0 & 1 & 0 & 0 & \cdots \\
0 & 0 & 1 & 0 & 0 & \cdots \\
& \vdots & & & \ddots \\
\end{pmatrix}}^q
\begin{pmatrix}
1 \\
0 \\
0 \\
0 \\
\vdots \\
\end{pmatrix}
}_{1}^{\log_3{2}}.
\end{equation*}
In Example \ref{example: calculation involving infinite matrix}, we will calculate the Hausdorff dimension of $X$ to be $\log_2{\beta} = 1.20\cdots$, where $\beta$ is the spectral radius of the infinite matrix
\begin{equation*}
\begin{pmatrix}
b_0  & b_1 &  b_2 & \cdots \\
1 & 0 & 0 & \cdots \\
0 & 1 & 0 & 0 & \cdots \\
0 & 0 & 1 & 0 & \cdots \\
& \vdots & & \ddots \\
\end{pmatrix}.
\end{equation*}
Or equivalently, $\beta$ is the unique positive solution to the equation
\[ \beta = \sum_{q = 0}^{\infty} \frac{b_q}{\beta^q} = 1 + \frac{1}{\beta} + \frac{2^{\log_4{3}}}{\beta^2} + \frac{2 + 2^{\log_4{3}}}{\beta^3} + \frac{3 + 2^{\log_4{3}} + 3^{\log_4{3}}}{\beta^4} + \cdots = 2.30\cdots. \]
The calculation is done by partitioning the summands by direction to express the Hausdorff dimension using the products of two infinite matrices. Then, we perform a similar trick again to denote the sum as a power of a single infinite matrix, revealing the connection with its spectral radius.

\section{Combinatorial expression for the Hausdorff dimension} \label{section: results}

In section \ref{section: proof}, we will prove the new formula for the Hausdorff dimension of self-affine fractals. Here, we will state this result and use it to investigate the dimension theory of some self-affine fractals.

\begin{notation}
Fix a natural number $r$. Consider a subshift $S \subset D_r^{\mathbb{N}}$, meaning $\sigma(S) \subset S$ for the shift map $\sigma(e^{(1)}, e^{(2)}, \cdots) = (e^{(2)}, e^{(3)}, \cdots)$. Let $m_1 \leq m_2 \leq \cdots \leq m_r$ be natural numbers, and $I_i = \{0, 1, \ldots, m_i-1\}$. Let $D_i = I_1 \times I_2 \times \cdots \times I_i$ for $1 \leq i \leq r$. For a natural number $N$, denote by $S|_N \subset D_r^N$ the restriction of $S$ to its first $N$ coordinates. Let $S^{(r)} = S$, and inductively define $p_i: S^{(i+1)}|_N \rightarrow D_i^N$ as the projection:
\begin{equation*}
p_i \left( \left( e^{(n)}_1, e^{(n)}_2, \ldots, e^{(n)}_{i}, e^{(n)}_{i+1} \right)_{n=1}^N \right) = \left( e^{(n)}_1, e^{(n)}_2, \ldots, e^{(n)}_{i} \right)_{n=1}^N,
\end{equation*}
We will abuse this notation for any $N \in \mathbb{N}$ or $N = \mathbb{N}$. Let $S^{(i)} = p_i(S^{(i+1)})$. Set $\hspace{2pt} D_i^N(s_{i-1}^{}) = \left\{ s \in D_i^N \setcond p_{i-1}^{}(s) = s_{i-1}^{} \right\} \hspace{2pt}$ for $\hspace{2pt} 2 \leq i \leq r \hspace{2pt}$ and $\hspace{2pt} s_{i-1}^{} \in D_{i-1}^N$. Also, define $R: D_r^{\mathbb{N}} \rightarrow {\mathbb{T}}^r$ by
\begin{equation*}
R \left( \left( e^{(n)}_1, e^{(n)}_2, \ldots, e^{(n)}_r \right)_{n=1}^{\infty} \right) = \left( \sum_{k=0}^{\infty} \frac{e^{(k)}_1}{{m_1}^k}, \cdots, \sum_{k=0}^{\infty} \frac{e^{(k)}_r}{{m_r}^k} \right).
\end{equation*}
Then, $R(S)$ is a compact set invariant under the map $T_A: {\mathbb{T}}^r \rightarrow {\mathbb{T}}^r$, which is defined as multiplying the matrix $A = \mathrm{diag}(m_1, m_2, \ldots, m_r)$. Conversely, if $K \subset {\mathbb{T}}^r$ is a $T_A$-invariant compact set, there is a subshift of $S \subset D_r^{\mathbb{N}}$ with $K = R(S)$, which is called the \textbf{symbolic model} of $K$. We will abuse $R$ to denote this map for any $r$ from now on.
\end{notation}

The following theorem is our main result, the proof of which is given in section \ref{section: proof}.

\begin{theorem} \label{theorem: combinatorial expression of Hausdorff dimension}
We have the following formula for any $T_A$-invariant compact set $K \subset {\mathbb{T}}^r$ for its Hausdorff dimension.
\begin{align} \nonumber
\mathrm{dim}_{\mathrm{H}}(X)
&= \lim_{N \to \infty} \spa \frac{1}{N} \spa \log_{m_1}{\hspace{-3pt}\sum_{s_1^{} \in S^{(1)}|_N} \hspace{-3pt}
\left( \sum_{s_2^{} \in p_1^{-1}(s_1^{})} \hspace{-5pt}
	\left( \cdots
		\left( \sum_{s_{r-1}^{} \in p_{r-2}^{-1}(s_{r-2}^{})} \hspace{-12pt}
		{\abs{p_{r-1}^{-1}(s_{r-1}^{})}}^{a_{r-1}^{}}
		\right)^{a_{r-2}^{}} \hspace{-9pt} \cdots
	\right)^{a_2^{}}
\right)^{a_1^{}}} \\[4pt] \nonumber
&= \lim_{N \to \infty} \spa \frac{1}{N} \spa \log_{m_1}{\hspace{-1pt}\sum_{s_1^{} \in D_1^N} \hspace{-1pt}
\left( \sum_{s_2^{} \in D_2^N(s_1^{})} \hspace{-3pt}
	\left( \cdots
		\left( \sum_{s_{r-1}^{} \in D_{r-1}^N(s_{r-2}^{})} \hspace{-10pt}
		{\abs{p_{r-1}^{-1}(s_{r-1}^{})}}^{a_{r-1}^{}}
		\right)^{a_{r-2}^{}} \hspace{-9pt} \cdots
	\right)^{\hspace{-2pt} a_2^{}}
\right)^{a_1^{}}}.
\end{align}
Here, $K = R(S)$, where $S \subset D_r^{\mathbb{N}}$ is the symbolic model of $K$, and $a_i^{} = \log_{m_{i+1}^{}}{m_i^{}} \hspace{2pt}$ for $ \spa 1 \leq i \leq r-1$.
\end{theorem}

\begin{breakremark}
\quad(1) Note that if $s_{r-1}^{} \in D_{r-1}^N(s_{r-2}^{})$ is ``illegal'' in the sense that $s_{r-1}^{} \notin S^{(r-1)}|_N$, then \\ $\abs{p_{r-1}^{-1}(s_{r-1}^{})} = 0$. Therefore, the second line is immediate from the first line.

\noindent \quad(2) The Hausdorff dimension formula of self-affine sponges in \cite[Theorem 1.2]{Kenyon--Peres} readily follows from this theorem. See the proof in \cite[Claim 1.6]{Alibabaei} for details.
\end{breakremark}

\begin{example} \label{example: uniform growth assumption}
Suppose that for every $2 \leq i \leq r$ and natural number $N$, there is a sequence $\left( P_N^{(i)} \right)_{N \in \mathbb{N}}$ such that we have a uniform convergence
\[ \lim_{N \to \infty} {\left( \frac{\abs{p_{i-1}^{-1}(s_{i-1}^{(N)})}}{P_N^{(i)}} \right)}^{\frac{1}{N}} = 1 \]
for any $s_{i-1}^{} \in S^{(i-1)}$, where $s_{i-1}^{(N)}$ is the first $N$ coordinates of $s_{i-1}^{}$ (note that $p_{i-1}^{-1}(s_{i-1}^{})$ is non-empty). We then say that $S$ has {\it{\textbf{uniformly growing word complexity}}}. Also, let $P^{(1)}_N = S^{(1)}|_N$. We have
\begin{flalign*}
& \mathrm{dim}_{\mathrm{H}}(X \hspace{-1pt}) \hspace{-3pt}
= \hspace{-2pt} \lim_{N \to \infty} \spa \frac{1}{N} \spa \log_{m_1^{}}{
\left(
	P_N^{(1)} {P_N^{(2)}}^{a_1} {P_N^{(3)}}^{a_1a_2} \cdots {P_N^{(r)}}^{a_1\cdots a_{r-1}}
\right)} &
\end{flalign*}
\begin{flalign*}
&= \lim_{N \to \infty} \frac{1}{N}
\left( \log_{m_1^{}}{ \hspace{-2pt}
			\left(
				P_N^{(1)} \hspace{-2pt} \cdots P_N^{(r)}
			\right)^{ \hspace{-2pt} a_1 \cdots a_{r-1}}}  \hspace{-4pt} +
			\sum_{i = 2}^{r-1}
			\log_{m_1^{}}{ \hspace{-2pt}
				\left(
					P_N^{(1)} \hspace{-2pt} \cdots P_N^{(i)}
				\right)^{ \hspace{-3pt} (1-a_i)a_1 \cdots a_{i-1}}} \hspace{-6pt} +
			\log_{m_1^{}}{ \hspace{-2pt}
				\left(
					P_N^{(1)}
				\right)^{ \hspace{-3pt} 1-a_1}} \hspace{-2pt}
\right) \\[5pt]
&= \lim_{N \to \infty} \spa
\left(\frac{1}{N \log{m_r^{}}} \spa \log{
	\left( \# S|_N 
	\right)} +
	\sum_{i = 1}^{r-1}
	\left( \frac{1}{\log{m_i^{}}} - \frac{1}{\log{m_{i+1}^{}}} \right)
	\frac{1}{N} \spa \log{\left( \# S^{(i)}|_N \right)}
\right) \\[1pt]
&= \frac{1}{\log{m_r}} h_{\mathrm{top}}(S, \sigma)
+ \sum_{i = 1}^{r-1} \left( \frac{1}{\log{m_i}} - \frac{1}{\log{m_{i+1}}} \right) h_{\mathrm{top}}(S^{(i)}, \sigma).
\end{flalign*}
\end{example}

Under this assumption, the corresponding self-affine fractal presumably has the same Minkowski and Hausdorff dimension. This is partially explained in Corollary \ref{corollary: uniform complexity implies no dimension hiatus} for sofic sets.

\subsection{Dimension theory of sofic sets} \label{section: dimension theory of sofic sets}

We will now focus on sofic sets. Weiss \cite{Weiss} defined {\it{sofic systems}} as subshifts which are factors of shifts of finite type. Using the results by \cite{Boyle--Kitchens--Marcus}, sofic systems can also be defined as follows.

\begin{definition}[{{\cite[Proposition 3.6]{Kenyon--Peres: sofic}}}] \label{definition: sofic systems}
Consider a finite directed graph $G = \langle V, E \rangle$ with possible loops and multiple edges. Let $I$ be a set of labels. Suppose that edges of $G$ are labeled using elements of $I$ in a ``right-resolving'' fashion: every two edges emanating from the same vertex have different labels. Then, the set $S \subset I^{\mathbb{N}}$ of sequences of labels that arise from infinite paths in $G$ is called the \textbf{sofic system}.
\end{definition}

For $I = D_r$, the image $R(S) \subset \mathbb{T}^r$ is called the \textbf{sofic set}. When $r=2$, Kenyon and Peres \cite{Kenyon--Peres: sofic} proved that the Hausdorff dimension of sofic sets can be written as the limit of a certain combinatorial sum of matrix products. As we shall explain, Theorem \ref{theorem: combinatorial expression of Hausdorff dimension} contains the direct generalization of this result to arbitrary $r$.


Suppose $S$ is a sofic system defined from a digraph with vertices $\{v_1, v_2, \ldots, v_n\}$. For $s \in D_{r-1}$, we define the restricted adjacency matrix $A_{s} = \left( a_{ij}(s) \right)_{i,j} \in \mathrm{M}_n(\mathbb{R})$ by
\[ a_{ij}(s) = \abs{ \left\{ 
				e \in E \setcond \text{$e$ is from vertex $v_i$ to $v_j$ and its label is $(s, k) \in D_r$ for some $k \in I_r$} 
			\right\}}. \]
Fix $s_{r-2}^{} \in S^{(r-2)}|_N$ , and take any $s_{r-1}^{} \in p_{r-2}^{-1}(s_{r-2}^{})$. If $s_{r-1} = (s^{(1)}, s^{(2)}, \ldots, s^{(N)})$, where $s^{(i)} \in D_{r-1}$ for all $i$. Then, by the right-resolving property,
\[ \abs{p_{r-1}^{-1}(s_{r-1}^{})} \asymp \norm{ A_{s^{(1)}} A_{s^{(2)}} \cdots A_{s^{(N)}} }. \]
Here, $A \asymp B $ means there is a constant $c > 0$ independent of $N$ with $c^{-1}B \leq A \leq cB$, and the norm of the matrix is the sum of all entries. We have the following corollary, which is the direct generalization of \cite[Theorem1.1]{Kenyon--Peres: sofic}.

\begin{corollary} \label{corollary: combinatorial formula for the Hausdorff dimension of sofic sets}
For any sofic system $S \subset D_r^{\mathbb{N}}$, we have
\begin{equation*} 
\mathrm{dim}_{\mathrm{H}}(X)
= \lim_{N \to \infty} \spa \frac{1}{N} \log_{m_1^{}}{\hspace{-5pt}\sum_{s_1^{} \in D_1^N} \hspace{-3pt}
\left( \hspace{-1pt} \sum_{s_2^{} \in D_2^N(s_1^{})} \hspace{-3pt}
	\left( \hspace{-2pt} \cdots \hspace{-2pt}
		\left( \sum_{\substack{
					s_{r-1}^{} \in D_{r-1}^N(s_{r-2}^{}) \\
					s_{r-1}^{} = (s^{(1)}, \ldots, s^{(N)})}} \hspace{-12pt}
		{\norm{ A_{s^{(1)}} \cdots A_{s^{(N)}} } }^{a_{r-1}^{}} \hspace{-2pt}
		\right)^{\hspace{-2pt} a_{r-2}^{}} \hspace{-14pt} \cdots
	\right)^{\hspace{-2pt} a_2^{}} \hspace{1pt}
\right)^{\hspace{-2pt} a_1^{}}}.
\end{equation*}
Here, the norm for matrices can be arbitrary.
\end{corollary}

We will later give some examples where we can explicitly calculate the limit above. Before that, let us focus on the theory of the Minkowski dimension (Box dimension) for sofic sets. Let $A$ be the adjacency matrix of $G$ in the usual sense, then $A = \sum_{s \in I_{r-1}} A_s$.
\begin{proposition} \label{proposition: Minkowski dimension}
Consider a sofic system $S$ and suppose its adjacency matrix $A$ is \textbf{primitive}, meaning some power of $A$ has positive entries. Then,  the sofic set $X = R(S)$ has the following Minkowski dimension.
\begin{align*}
\mathrm{dim}_{\mathrm{M}}(X)
&=
\frac{1}{\log{m_r}} h_{\mathrm{top}}(S, \sigma)
+ \sum_{i = 1}^{r-1} \left( \frac{1}{\log{m_i}} - \frac{1}{\log{m_{i+1}}} \right) h_{\mathrm{top}}(S^{(i)}, \sigma),
\end{align*}
where $\sigma$ is the left-shift map. We also have $h_{\mathrm{top}}(S, \sigma) = \log{\rho(A)}$ with $\rho(A)$ being the spectral radius of $A$.
\end{proposition}

\begin{proof}
The method in \cite[Proposition 3.5]{Kenyon--Peres: sofic} for the $2$-dimensional case works here as well. We include proof here for the sake of completeness. Let $L_r: E \rightarrow D_r$ be the labelling of $G$. Define $\phi_i: D_r \rightarrow D_i$ by
\[ \phi_i (e_1, \ldots, e_r) = (e_1, \ldots, e_i), \]
and let $L_i: E \rightarrow D_i$ be the composition $\phi_i \circ L_r$ for $1 \leq i \leq r$. Define
\[ N_i(k) = \# \left\{ \left( L_i(e_1), \ldots, L_i(e_k) \right) \setcond \text{$(e_1, \ldots, e_k)$ is a path in $G$} \right\} = \# \left( S^{(i)}|_k \right). \]

Let $b_i = \log_{m_i} {m_1}$ and $M(k)$ be the number of distinct sequences of the form
\[ \Big( L_r(e_1), \ldots, L_r(e_{\lfloor b_rk \rfloor}), \hspace{4pt}
				L_{r-1}(e_{\lfloor b_rk \rfloor + 1}), \ldots, L_{r-1}(e_{\lfloor b_{r-1}k \rfloor}), \hspace{4pt} \ldots \hspace{4pt},
				L_1(e_{\lfloor b_2k \rfloor} + 1), \ldots, L_1(e_k)
			 \Big), \]
where $(e_1, \ldots, e_k)$ is a path in $G$. We have
\begin{equation} \label{equation: topological entropy of factors}
h_{\mathrm{top}}(S^{(i)}) = \lim_{k \to \infty} \frac{1}{k} \log{ \left( N_i(k) \right) } \quad (\text{\it{ for }} 1 \leq i \leq r)
\end{equation}
and
\begin{equation} \label{equation: Minkowski dimension of X in the proof}
\mathrm{dim}_{\mathrm{M}}(X) = \lim_{k \to \infty} \frac{1}{k} \log_{m_1} {M(k)}.
\end{equation}

The formulae (\ref{equation: topological entropy of factors}) are immediate from the definition of topological entropy. The equation (\ref{equation: Minkowski dimension of X in the proof}) follows since a cube of size $m_1^k \times m_2^{\lfloor b_2k \rfloor} \times \cdots \times m_r^{\lfloor b_rk \rfloor}$ is approximately a cube with side length $m_1^{}$, and $M(k)$ is the number of such cubes intersecting $X$.

By definition,
\[ M(k) \leq N_r\left( \lfloor b_r k \rfloor \right)
\cdot N_{r-1}\left( \lfloor b_{r-1}k \rfloor - \lfloor b_r k \rfloor \right) \hspace{3pt} \cdots \hspace{3pt}
N_1\left( k - \lfloor b_2 k \rfloor \right). \]
The other direction is due to the assumption that $A$ is primitive; suppose $A^d$ is positive, then,
\[ M(k) \geq N_r\left( \lfloor b_r k \rfloor \right) \cdot
N_{r-1}\left( \lfloor b_{r-1}k \rfloor - \lfloor b_r k \rfloor - d \right) \hspace{3pt} \cdots \hspace{3pt}
N_1\left( k - \lfloor b_2 k \rfloor - d \right). \]

Combining these inequalities gives us
\begin{flalign*}
&
\begin{multlined}[t][11.5cm]
\frac{1}{k} \log{N_r\left( \lfloor b_r k \rfloor \right)} + \frac{1}{k} \log{N_{r-1}\left( \lfloor b_{r-1}k \rfloor - \lfloor b_r k \rfloor - d \right)} + \cdots +
\frac{1}{k} \log{N_1\left( k - \lfloor b_2 k \rfloor - d \right)} \\
\leq M(k) \leq 
\frac{1}{k} \log{N_r\left( \lfloor b_r k \rfloor \right)} + \frac{1}{k} \log{N_{r-1}\left( \lfloor b_{r-1}k \rfloor - \lfloor b_r k \rfloor \right)} + \cdots +
\frac{1}{k} \log{N_1\left( k - \lfloor b_2 k \rfloor \right)}.
\end{multlined} &
\end{flalign*}
Letting $k \to \infty$ and using formulae (\ref{equation: topological entropy of factors}) and (\ref{equation: Minkowski dimension of X in the proof}) yields the desired equation.

Also, by the right-resolving property, we have $\frac{1}{|V|} \norm{A^k} \leq N_r(k) \leq \norm{A^k}$. This implies \\ $h_{\mathrm{top}}(S, \sigma) = \log{\rho(A)}$. (Note that the topological entropy of symbolic spaces $S^{(i)}$ may not be the spectral radius of the adjacency matrix of the corresponding digraph. This is because the right-resolving property is usually lost after projections.)
\end{proof}

The assumption that $A$ is primitive cannot be dropped: see the remark after \cite[Proposition 3.5]{Kenyon--Peres: sofic} for details.

This proposition, together with the calculation in Example \ref{example: uniform growth assumption}, suggests that the uniform growing assumption for word complexity is the defining condition for sofic systems to have the same Hausdorff and Minkowski dimension. For sofic sets in $\mathbb{R}^2$, the assumption that word complexity grows uniformly seems to correspond to the following definition.
\begin{definition}
A sofic system $S \subset D_2^{\mathbb{N}}$ is said to have {\it{\textbf{uniform complexity}}} if all $A_{s}$ for $s \in S^{(r-1)}|_1$ have the common Perron-Frobenius eigenvalue and the corresponding eigenvector; there is a positive number $\lambda$ and a column (or row) vector $v$ with positive entries such that $v \spa A_{s} = \lambda v$ (or $A_{s} v = \lambda v$) for every $s \in S^{(r-1)}|_1$.
\end{definition}

By the following Corollary \ref{corollary: uniform complexity implies no dimension hiatus}, if $S$ has uniform complexity, then $X = R(S)$ has the same Minkowski and Hausdorff dimension.

\begin{figure}[h!]
\includegraphics[width=15cm]{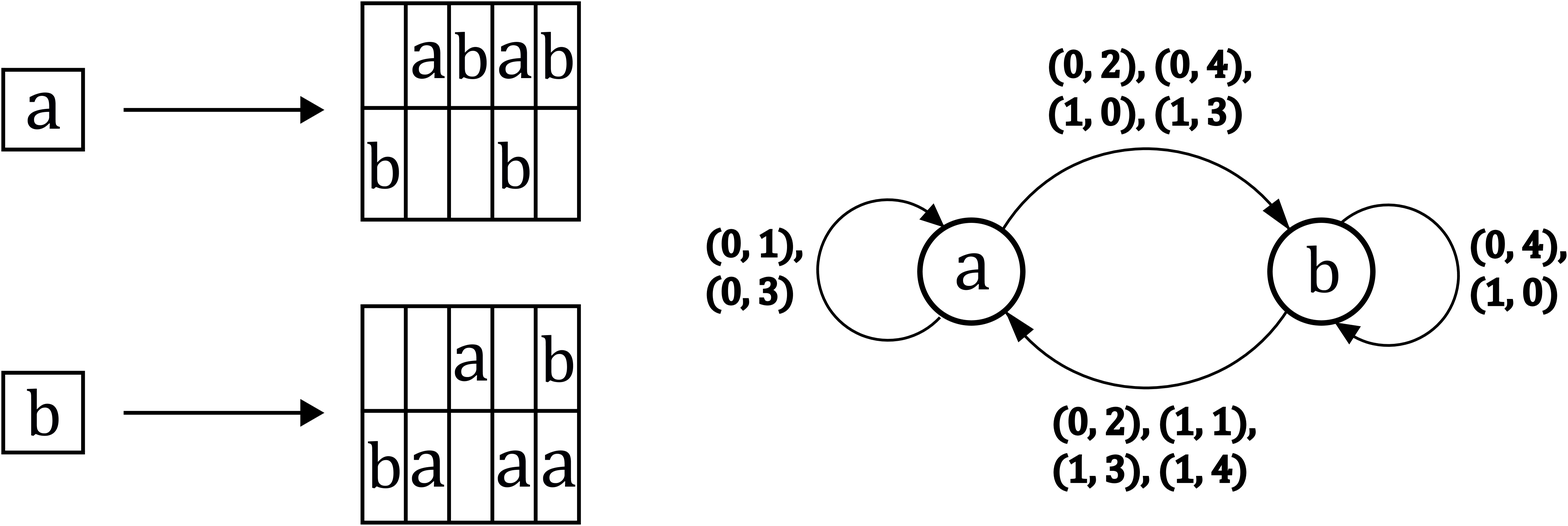}
\caption{A classic example of a sofic set with no dimension hiatus} \label{figure: coincidence of definition 2}
\end{figure}

All the examples of sofic sets considered in the work of Olivier \cite{Olivier} and Kenyon-Peres \cite{Kenyon--Peres: sofic} with no dimension hiatus satisfy the condition above. In fact, they all satisfy the following stronger condition: the number of each symbol in each row in the initial pattern is equal for all rows. In this case, all matrices have $(1, 1, \ldots, 1)$ as a shared positive eigenvector with the same eigenvalue, implying uniform complexity. Figure \ref{figure: coincidence of definition 2} is such an example, with
\begin{equation*}
A_0 =
\begin{pmatrix}
2 & 2 \\
1 & 1 \\
\end{pmatrix},
\hspace{4pt} A_1 =
\begin{pmatrix}
0 & 2 \\
3 & 1 \\
\end{pmatrix}.
\end{equation*}
We immediately calculate its dimension to be $1 + \log_5{3}$.

Figure \ref{figure: coincidence of definition} in the introduction also satisfies uniform complexity, where $(2, 1)$ is the common Perron-Frobenius eigenvector with the same eigenvalue for the adjacency matrices
\begin{equation*}
A_0 =
\begin{pmatrix}
2 & 1 \\
0 & 0 \\
\end{pmatrix},
\hspace{4pt} A_1 =
\begin{pmatrix}
0 & 1 \\
4 & 0 \\
\end{pmatrix}.
\end{equation*}

\begin{corollary} \label{corollary: uniform complexity implies no dimension hiatus}
Let $S \subset D_2^{\mathbb{N}}$ be a sofic system with primitive adjacency matrix. If $S$ has uniform complexity, then
\begin{equation*}
\mathrm{dim}_{\mathrm{H}}(X) = \mathrm{dim}_{\mathrm{M}}(X).
\end{equation*}
\end{corollary}

\begin{proof}
Suppose that a positive vector $v$ is the common left eigenvector with $\lambda$ as the eigenvalue for all $A_{s}$ for $s \in S^{(r-1)}|_1$. There are positive constants $c$ and $C$ with
\[ c (1, 1, \ldots, 1) \leq v \leq C (1, 1, \ldots, 1). \]
Then, for any $s =( s^{(i)} )_{i = 1}^{\infty} \in S^{(1)}$ and a natural number $N$, we have
\begin{align*}
\abs{p_1^{-1}(s^{(1)}, \ldots, s^{(N)})}
& \asymp \norm{ A_{s^{(1)}} \cdots A_{s^{(N)}}} \\
& \asymp \norm{ v A_{s^{(1)}} \cdots A_{s^{(N)}}} \\
& \asymp \lambda^N.
\end{align*}
Therefore, $S$ has uniformly growing word complexity (see Example \ref{example: uniform growth assumption}). By the calculation in Example \ref{example: uniform growth assumption}, we have
\[ \mathrm{dim}_{\mathrm{H}}(X) = \frac{1}{\log{m_2}} h_{\mathrm{top}}(S, \sigma)
+ \left( \frac{1}{\log{m_1}} - \frac{1}{\log{m_2}} \right) h_{\mathrm{top}}(S^{(1)}, \sigma). \]
This is also the Minkowski dimension of $X$ by Proposition \ref{proposition: Minkowski dimension}.
\end{proof}

We conjecture that the converse is also true.
\begin{conjecture}
A sofic system $S \subset D_2^{\mathbb{N}}$ with primitive adjacency matrix has the same Hausdorff and Minkowski dimension if and only if it has uniform complexity.
\end{conjecture}

Let us consider a sofic system $S \subset D_r^{\mathbb{N}}$ with $r \geq 3$. Unlike planar cases, we cannot give an intuitive condition for $S$ to have uniformly growing word complexity; the projected graph may not be right-resolving, making it impossible to estimate $\abs{ p_1^{-1} }$ using matrices. \\[3.5pt]

We now calculate the Hausdorff dimension of some sofic sets in $\mathbb{R}^3$ or higher using Corollary \ref{corollary: combinatorial formula for the Hausdorff dimension of sofic sets}. We can calculate the Hausdorff dimension if all $A_s$ for $s \in D_{r-1}$ have a shared positive eigenvector.

\begin{proposition}
Let $S$ be a sofic system, and suppose that there is a positive vector $v$ with $A_s \spa v = \lambda_s \spa v$ for all $s \in D_{r-1}$. Then,
\begin{equation*}
\mathrm{dim}_{\mathrm{H}}(X) = \log_{m_1^{}}{\sum_{s_1^{} \in D_1} 
\left(  \sum_{s_2^{} \in D_2(s_1^{})} 
	\left(  \cdots
		\left( \sum_{s_{r-1}^{} \in D_{r-1}(s_{r-2}^{})} 
			\lambda_{s_{r-1}^{}}^{a_{r-1}}
		\right)^{a_{r-2}^{}} \cdots
	\right)^{ a_2^{}} 
\right)^{a_1^{}}}.
\end{equation*}
\end{proposition}

\begin{proof}
First, observe that the ratio of $\norm{A_{s^{(1)}} \cdots A_{s^{(N)}}}$ and 
\[ \norm{A_{s^{(1)}} \cdots A_{s^{(N)}} \spa v } = \lambda_{s^{(1)}} \cdots \lambda_{s^{(N)}} \spa \norm{v} \]
is bounded. We then have
\begin{flalign*}
&
\hspace{15pt} \mathrm{dim}_{\mathrm{H}}(X)
=
\lim_{N \to \infty} \spa \frac{1}{N} \log_{m_1^{}}{\hspace{-5pt}\sum_{s_1^{} \in D_1^N} 
\left( \cdots 
		\left( \sum_{\substack{
					s_{r-1}^{} \in p_{r-2}^{-1}(s_{r-2}^{}) \\
					s_{r-1}^{} = (s^{(1)}, \ldots, s^{(N)})}} \hspace{-8pt}
		{\norm{A_{s^{(1)}} \cdots A_{s^{(N)}} \spa v}}^{a_{r-1}^{}} \hspace{-2pt}
		\right)^{\hspace{-2pt} a_{r-2}^{}} \hspace{-9pt} \cdots \hspace{1pt}
\right)^{\hspace{-2pt} a_1^{}}} &
\end{flalign*}
\begin{flalign*}
& \hspace{58.5pt} =
\lim_{N \to \infty} \spa \frac{1}{N} \log_{m_1^{}}{\hspace{-5pt}\sum_{s_1^{} \in D_1^N} 
\left( \cdots 
		\left( \sum_{\substack{
					s_{r-1}^{} \in p_{r-2}^{-1}(s_{r-2}^{}) \\
					s_{r-1}^{} = (s^{(1)}, \ldots, s^{(N)})}} \hspace{-8pt}
		\lambda_{s^{(1)}}^{a_{r-1}^{}} \cdots \lambda_{s^{(N)}}^{a_{r-1}^{}} \hspace{-2pt}
		\right)^{\hspace{-2pt} a_{r-2}^{}} \hspace{-9pt} \cdots \hspace{1pt}
\right)^{\hspace{-2pt} a_1^{}}} &
\end{flalign*}
\begin{flalign*}
& \hspace{58.5pt} =
\lim_{N \to \infty} \spa \frac{1}{N} \log_{m_1^{}}{
{\left\{ \sum_{s_1^{} \in D_1} 
\left(  \sum_{s_2^{} \in D_2(s_1^{})} 
	\left(  \cdots
		\left( \sum_{s_{r-1}^{} \in D_{r-1}(s_{r-2}^{})} 
			\lambda_{s_{r-1}^{}}^{a_{r-1}}
		\right)^{a_{r-2}^{}} \cdots
	\right)^{ a_2^{}} 
\right)^{a_1^{}}
\right\}}}^N &
\end{flalign*}
\begin{flalign*}
& \hspace{58.5pt} = \log_{m_1^{}}{\sum_{s_1^{} \in D_1} 
\left(  \sum_{s_2^{} \in D_2(s_1^{})} 
	\left(  \cdots
		\left( \sum_{s_{r-1}^{} \in D_{r-1}(s_{r-2}^{})} 
			\lambda_{s_{r-1}^{}}^{a_{r-1}}
		\right)^{a_{r-2}^{}} \cdots
	\right)^{ a_2^{}} 
\right)^{a_1^{}}}. &
\end{flalign*}
\end{proof}

Following the idea of Kenyon and Peres in \cite{Kenyon--Peres: sofic}, we calculate the Hausdorff dimension of some sofic sets by partitioning the summands by direction. The following is an example with seemingly simple matrices, demonstrating how complex the calculation for the Hausdorff dimension can be.

\begin{example} \label{example: calculation involving infinite matrix}
Consider the sofic system defined from the digraph in Figure \ref{figure: digraph for sofic set with calculation using infinite matrix}. We have
\begin{equation*}
A_{(0, 1)} =
\begin{pmatrix}
0 & 0 \\
0 & 1 \\
\end{pmatrix},
\hspace{4pt} A_{(1, 0)} =
\begin{pmatrix}
0 & 1 \\
0 & 0 \\
\end{pmatrix},
\hspace{4pt} A_{(1, 1)} =
\begin{pmatrix}
1 & 0 \\
1 & 1 \\
\end{pmatrix}.
\end{equation*}

First, note that $A_{(1, 1)} A_{(1, 0)} \begin{psmallmatrix} 0 \\ 1 \end{psmallmatrix} = \begin{psmallmatrix} 1 \\ 1 \end{psmallmatrix}$. This implies that (for $a_1^{} = \log_3{2}$ and $a_2^{} = \log_4{3}$)
\begin{align*}
\sum_{s_1^{} \in \{0, 1\}^N} 
\left( \spa \sum_{\substack{
					s_{2}^{} \in D_{2}^N(s_{1}^{}) \\
					s_{2}^{} = (s^{(1)}, \ldots, s^{(N)})}} \hspace{-2pt}
		{\left[ (1,1)A_{s^{(1)}} \cdots A_{s^{(N)}}
			\begin{pmatrix}
			0 \\
			1 \\
			\end{pmatrix}
		\right] }^{a_{2}^{}}
\right)^{\hspace{-2pt} a_1^{}}
\end{align*}
\begin{align*}
\hspace{5pt} \leq \sum_{s_1^{} \in \{0, 1\}^N} 
\left( \spa \sum_{\substack{
					s_{2}^{} \in D_{2}^N(s_{1}^{}) \\
					s_{2}^{} = (s^{(1)}, \ldots, s^{(N)})}} \hspace{-2pt}
		{\left[ (1,1)A_{s^{(1)}} \cdots A_{s^{(N)}}
			\begin{pmatrix}
			1 \\
			1 \\
			\end{pmatrix}
		\right] }^{a_{2}^{}}
\right)^{\hspace{-2pt} a_1^{}}
\end{align*}
\begin{flalign*}
\hspace{58pt} = \sum_{s_1^{} \in \{0, 1\}^N} 
\left( \spa \sum_{\substack{
					s_{2}^{} \in D_{2}^N(s_{1}^{}) \\
					s_{2}^{} = (s^{(1)}, \ldots, s^{(N)})}} \hspace{-2pt}
		{\left[ (1,1)A_{s^{(1)}} \cdots A_{s^{(N)}} A_{(1,1)} A_{(1,0)}
			\begin{pmatrix}
			0 \\
			1 \\
			\end{pmatrix}
		\right] }^{a_{2}^{}}
\right)^{\hspace{-2pt} a_1^{}}
\end{flalign*}
\begin{flalign*}
\hspace{40pt} \leq \sum_{s_1^{} \in \{0, 1\}^{N+2}} 
\left( \spa \sum_{\substack{
					s_{2}^{} \in D_{2}^{N+2}(s_{1}^{}) \\
					s_{2}^{} = (s^{(1)}, \ldots, s^{(N+2)})}} \hspace{-2pt}
		{\left[ (1,1)A_{s^{(1)}} \cdots A_{s^{(N+2)}}
			\begin{pmatrix}
			0 \\
			1 \\
			\end{pmatrix}
		\right] }^{a_{2}^{}}
\right)^{\hspace{-2pt} a_1^{}}.
\end{flalign*}
Therefore, we have
\begin{equation*}
\mathrm{dim}_{\mathrm{H}}(X)
=
\lim_{N \to \infty} \spa \frac{1}{N} \log_{2}{\hspace{-3pt}\sum_{s_1^{} \in \{0, 1\}^N} 
\left( \spa \sum_{\substack{
					s_{2}^{} \in D_{2}^N(s_{1}^{}) \\
					s_{2}^{} = (s^{(1)}, \ldots, s^{(N)})}} \hspace{-2pt}
		{\left[ (1,1)A_{s^{(1)}} \cdots A_{s^{(N)}}
			\begin{pmatrix}
			0 \\
			1 \\
			\end{pmatrix}
		\right] }^{a_{2}^{}}
\right)^{\hspace{-2pt} a_1^{}}}.
\end{equation*}

Now, we partition the summands by direction. Namely, fix $s_1^{} \in \{0, 1\}^N$ and observe that for any $N$-tuple $(s^{(1)}, \ldots, s^{(N)}) \in D_2^N(s_1^{})$, the vector $(1, 1) A_{s^{(1)}} \cdots A_{s^{(N)}}$ is some constant times $(q, 1)$, where $q$ is a non-negative integer. Let $H_N(q, s_1^{})$ be the set of $(s^{(1)}, \ldots, s^{(N)})$ with such $q$. Let
\begin{equation*}
\Phi_N(q, s_1^{}) = \sum_{(s^{(1)}, \ldots, s^{(N)}) \in H_N(q, s_1^{})} {\left[ (1,1)A_{s^{(1)}} \cdots A_{s^{(N)}}
			\begin{pmatrix}
			0 \\
			1 \\
			\end{pmatrix}
		\right] }^{a_{2}^{}}.
\end{equation*}
Also, let $\Phi_0 = (1, 0, 0, \cdots)^{\mathrm{T}}$. Consider $(s_1^{}, 1) \in \{0, 1\}^{N+1}$. Since $(q, 1) A_{(1, 0)} = (0, q) = q (0, 1)$, and $(q, 1) A_{(1, 1)} = (q+1, 1)$, we see that
\begin{gather*}
\Phi_{N+1}\left( q+1 , (s_1^{}, 1) \right) = \Phi_N (q, s_1^{} ), \\
\Phi_{N+1}\left( 0 , (s_1^{}, 1) \right) = \sum_{q = 0}^{\infty} q^{a_2} \spa \Phi_N (q, s_1^{} ).
\end{gather*}
Or equivalently, letting
\begin{equation*}
M_1 = 
\begin{pmatrix}
0 & 1 & 2^{a_2^{}} & 3^{a_2^{}} & 4^{a_2^{}} & \cdots \\
1 & 0 & 0 & \cdots \\
0 & 1 & 0 & 0 & \cdots \\
0 & 0 & 1 & 0 & 0 & \cdots \\
& \vdots & & & \ddots \\
\end{pmatrix}
\end{equation*}
yields the following vector notation:
\[ \Phi_{N+1}\left(s_1^{}, 1\right) = M_1 \Phi_{N}\left(s_1^{}\right). \]

Similar argument works for $(s_1^{}, 0) \in \{0, 1\}^{N+1}$: since $(q, 1) A_{(0, 1)} = (0, 1)$, letting
\begin{equation*}
M_0 = 
\begin{pmatrix}
1 & 1 & 1 & \cdots \\
0 & 0 & 0 & \cdots \\
0 & 0 & 0 & \cdots \\
& \vdots & & \ddots \\
\end{pmatrix},
\end{equation*}
we have
\[ \Phi_{N+1}\left(s_1^{}, 0\right) = M_0 \Phi_{N}\left(s_1^{}\right). \]

This gives us the following expression for the Hausdorff dimension.
\begin{equation*}
\mathrm{dim}_{\mathrm{H}}(X)
=
\lim_{N \to \infty} \spa \frac{1}{N} \log_{2}{\hspace{-3pt}\sum_{(u_1^{}, \ldots, u_N^{}) \in \{0, 1\}^N}
{\norm{ M_{u_N^{}} \spa M_{u_{N-1}^{}} \cdots \spa M_{u_1^{}} \spa \Phi_0 }
}^{a_1^{}}}.
\end{equation*}
Here, the norm is the $l^1$-norm. Let
\[v = (1, \spa 1, \spa 0, \hspace{2pt} -2^{a_2^{}}, \hspace{2pt} -2^{a_2^{}}-3^{a_2^{}}, \hspace{2pt} -2^{a_2^{}}-3^{a_2^{}}-4^{a_2^{}}, \cdots). \]
Then, $v \spa M_1 = v$ and $v \spa M_0 = (1, 1, 1, \cdots)$. Since we have
\[ v \spa \Phi_0 = 1 \geq 0, \]
it follows that
\begin{align*}
& \sum_{(u_1^{}, \ldots, u_N^{}) \in D_{1}^N}
{ \left[ v \spa M_{u_N^{}} \spa M_{u_{N-1}^{}} \cdots M_{u_1^{}} \spa \Phi_0
\right] }^{a_1^{}} \\
& \leq \sum_{(u_1^{}, \ldots, u_N^{}) \in D_{1}^N}
{ \left[ (1, 1, 1, \cdots) \spa M_{u_N^{}} \spa M_{u_{N-1}^{}} \cdots M_{u_1^{}} \spa \Phi_0
\right] }^{a_1^{}} \\
& = \sum_{(u_1^{}, \ldots, u_N^{}) \in D_{1}^N}
{ \left[ v \spa M_0^{} M_{u_N^{}} \spa M_{u_{N-1}^{}} \cdots M_{u_1^{}} \spa \Phi_0
\right] }^{a_1^{}} \\
& \leq \sum_{(u_1^{}, \ldots, u_{N+1}^{}) \in D_{1}^{N+1}}
{ \left[ v \spa M_{u_{N+1}^{}} \spa M_{u_{N}^{}} \cdots M_{u_1^{}} \spa \Phi_0
\right] }^{a_1^{}}.
\end{align*}
Thus, we obtain
\begin{equation*}
\mathrm{dim}_{\mathrm{H}}(X)
=
\lim_{N \to \infty} \spa \frac{1}{N} \log_{2}{\hspace{-3pt}\sum_{(u_1^{}, \ldots, u_N^{}) \in D_{1}^N}
{ \left[ v \spa M_{u_N^{}} \spa M_{u_{N-1}^{}} \cdots M_{u_1^{}} \spa \Phi_0
\right] }^{a_1^{}}}.
\end{equation*}
For any $(u_1^{}, \ldots, u_N^{}) \in D_{1}^N$, the vector $M_{u_N^{}} M_{u_{N-1}^{}} \cdots M_{u_1^{}} \Phi_0$ is some constant times $M_1^q \Phi_0$. Denote by $U_N(q)$ the collection of such $(u_1^{}, \ldots, u_N^{})$. Let
\[ \Psi_N(q) = \sum_{(u_1^{}, \ldots, u_N^{}) \in U_N(q)} { \left[ v \spa M_{u_N^{}} \spa M_{u_{N-1}^{}} \cdots M_{u_1^{}} \spa \Phi_0
\right] }^{a_1^{}}. \]
Note that we have
\[ M_0^{} \spa M_1^q \spa \Phi_0 = \norm{M_1^q \spa \Phi_0} \spa \Phi_0. \]
This implies the following relations with $b_q = {\norm{M_1^q \spa \Phi_0}}^{a_1^{}}$;
\begin{gather*}
\Psi_{N+1}(0) = \sum_{q = 0}^{\infty} b_q \spa \Psi_N(q), \\
\Psi_{N+1}(q+1) = \Psi_{N}(q).
\end{gather*}
As before, letting
\begin{equation*}
L = 
\begin{pmatrix}
b_0  & b_1 &  b_2 & \cdots \\
1 & 0 & 0 & \cdots \\
0 & 1 & 0 & 0 & \cdots \\
0 & 0 & 1 & 0 & \cdots \\
& \vdots & & \ddots \\
\end{pmatrix}
\end{equation*}
yields the vector notation
\[ \Psi_{N+1} = L \spa \Psi_{N}. \]
This gives rise to
\[ \mathrm{dim}_{\mathrm{H}}(X)
=
\lim_{N \to \infty} \spa \frac{1}{N} \log_{2}{
\norm{ L^N \Psi_0 }} \]
with $\Psi_0 = (1, 0, 0, \cdots)^{\mathrm{T}}$.

Define $\beta$ to be the unique positive solution to the equation
\[ \beta = \sum_{q = 0}^{\infty} \frac{b_q}{\beta^q}. \]
It is immediate that $u = (1, \beta^{-1}, \beta^{-2}, \cdots)^{\mathrm{T}}$ satisfies $L \spa u = \beta u$. This implies
\[ \norm{ L^N \spa \Psi_0 } \leq \norm{ L^N \spa u } = \beta^N u. \]
Conversely, applying the Perron-Frobenius theorem to the upper-left $m \times m$ submatrix $L_m$ of $L$ yields
\[ \liminf_{N \to \infty} \norm{ L^N \spa \Psi_0 }^{\frac{1}{N}}
\geq \liminf_{N \to \infty} \norm{ L_m^N \spa \Psi_0^{(m)} }^{\frac{1}{N}} = \beta_m. \]
Here, $\beta_m$ is the spectral radius of $L_m$ which satisfies
\[ \beta_m = \sum_{q = 0}^{m-1} \frac{b_q}{\beta_m^q}. \]
Since $\beta_m \to \beta$ when $m \to \infty$, we obtain
\[ \lim_{N \to \infty} {\norm{ L^N \spa \Psi_0 }}^{\frac{1}{N}} = \beta. \]

Therefore, we have
\[ \mathrm{dim}_{\mathrm{H}}(X) = \log_{2}{ \beta } = 1.20\cdots. \]
\end{example}

In this example, we had to do the ``partition trick'' twice. The following example requires us to do it only once, making it a straight extension of the calculation done by Kenyon and Peres \cite{Kenyon--Peres: sofic}.

\begin{example}
Let $D_3 = \{0, 1, 2\} \times \{0, 1, 2, 3\} \times \{0, 1, 2, 3, 4\}$, and consider the sofic system defined from the digraph in Figure \ref{figure: calculation for a sofic set in R^3 number 2}. 
\begin{figure}[h!]
\includegraphics[width=10cm]{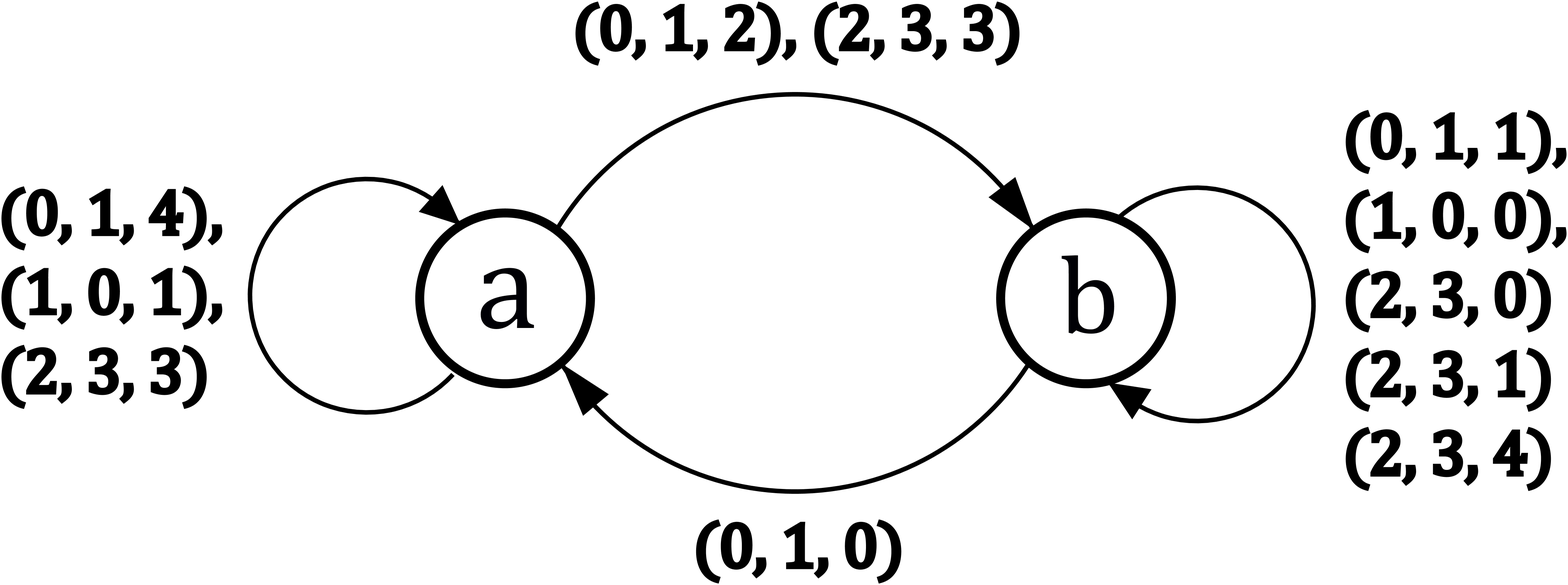}
\caption{A digraph for a sofic set in $\mathbb{R}^3$} \label{figure: calculation for a sofic set in R^3 number 2}
\end{figure}

\noindent We have
\begin{equation*}
A_{(0, 1)} =
\begin{pmatrix}
1 & 1 \\
1 & 1 \\
\end{pmatrix},
\hspace{4pt} A_{(1, 0)} =
\begin{pmatrix}
1 & 0 \\
0 & 1 \\
\end{pmatrix},
\hspace{4pt} A_{(2, 3)} =
\begin{pmatrix}
1 & 1 \\
0 & 3 \\
\end{pmatrix}.
\end{equation*}

Let $A_0 = A_{(0, 1)}$, $A_1 = A_{(1, 0)}$, and $A_2 = A_{(2, 3)}$. We have
\begin{equation*}
\mathrm{dim}_{\mathrm{H}}(X)
=
\lim_{N \to \infty} \spa \frac{1}{N} \log_{2}{ \sum_{\substack{
				s \in \{0, 1, 2\}^N \\
				s = (s^{(1)}, \ldots, s^{(N)})}}
			{\norm{A_{s^{(1)}} \cdots A_{s^{(N)}}}
			}^{a_{2}^{}a_1^{}}}.
\end{equation*}
Let $\alpha = a_1^{} a_2^{} = \log_5{3}$. We see that the vector $(1, 1) A_{s^{(1)}} \cdots A_{s^{(N)}}$ is some constant times $(1, \frac{3^q - 1}{2})$, where $q \geq 1$ is a natural number. Let $U^{(q)}_N \subset \{0, 1, 2\}^{N}$ be the collection of such $(s^{(1)}, \ldots, s^{(N)})$. Let
\begin{equation*}
\Phi_N(q) = \sum_{(s^{(1)}, \ldots, s^{(N)}) \in U^{(q)}_N} {\left[ (1,1)A_{s^{(1)}} \cdots A_{s^{(N)}}
			\begin{pmatrix}
			1 \\
			1 \\
			\end{pmatrix}
		\right] }^{\alpha}.
\end{equation*}
Also, let $\Phi_0 = (1, 0, 0, \cdots)^{\mathrm{T}}$. Since $(1, \frac{3^q - 1}{2}) A_0 = \frac{3^q + 1}{2} (1, 1)$, $(1, \frac{3^q - 1}{2}) A_1 = (1, \frac{3^q - 1}{2})$, and $(1, \frac{3^q - 1}{2}) A_2 = (1, \frac{3^{q+1} - 1}{2})$, we see that
\begin{gather*}
\Phi_{N+1}(q+1) = \Phi_N (q) + \Phi_N (q+1), \\
\Phi_{N+1}(1) = \sum_{q = 0}^{\infty} \left( \frac{3^q + 1}{2} \right)^{\alpha} \spa \Phi_N (q).
\end{gather*}
Letting
\begin{equation*}
M = 
\begin{pmatrix}
1 + 2^\alpha & 5^\alpha & 14^\alpha & \cdots & \left( \frac{3^q + 1}{2} \right)^{\alpha} & \cdots \\
1 & 1 & 0 & \cdots \\
0 & 1 & 1 & 0 & \cdots \\
0 & 0 & 1 & 1 & \\
& \vdots & & \ddots &\ddots \\
\end{pmatrix}
\end{equation*}
yields the following vector notation:
\[ \Phi_{N+1} = M \Phi_{N}. \]
We conclude that
\[ \mathrm{dim}_{\mathrm{H}}(X)
=
\lim_{N \to \infty} \spa \frac{1}{N} \log_{3}{
\norm{ M^N \Phi_0 }}. \]
Let $b_n$ be the $(1, n)$-th entry of $M$, and $\beta$ the unique positive solution of the equation
\[ \beta - \frac{1}{\beta} = b_1 + \sum_{k = 1}^{\infty} (b_{k+1} - b_k) \frac{1}{\beta^k}. \]
We then have
\begin{equation*}
M
\begin{pmatrix}
1 - \beta^{-1} \\
\beta^{-1} - \beta^{-2} \\
\beta^{-2} - \beta^{-3} \\
\vdots
\end{pmatrix}
=
\begin{pmatrix}
\beta - \beta^{-1} \\
1 - \beta^{-2} \\
\beta^{-1} - \beta^{-3} \\
\vdots
\end{pmatrix}
= (\beta + 1)
\begin{pmatrix}
1 - \beta^{-1} \\
\beta^{-1} - \beta^{-2} \\
\beta^{-2} - \beta^{-3} \\
\vdots
\end{pmatrix}.
\end{equation*}
This implies that $\rho(M) = \beta + 1$. By the same argument as before, it follows that
\[ \mathrm{dim}_{\mathrm{H}}(X) = \log_{3}{ \left( \beta + 1 \right)} = 1.231\cdots. \]

\end{example}

\section{Proof} \label{section: proof}

\subsection{A brief introduction to weighted topological entropy}

In this section, we will prove the main theorem of this paper. For that, we will be using weighted topological entropy, a concept recently introduced by Feng and Huang \cite{Feng--Huang}. Following the work of Tsukamoto \cite{Tsukamoto}, who provided a new definition of it in the simplest case, we gave rise to a new definition of weighted topological entropy under general settings in \cite{Alibabaei}. In this paper, we will be using both definitions of this invariant. Let us briefly review the notations, starting from the basics. Refer to the book of Walters \cite{Walters} for the details.

A pair $(X, T)$ is called a \textbf{dynamical system} if $X$ is a compact metrizable space and $T: X \rightarrow X$ is a continuous map. A map $\pi: X \rightarrow Y$ between dynamical systems $(X, T)$ and $(Y, S)$ is said to be a \textbf{factor map} if $\pi$ is a continuous surjection and $\pi \circ T = S \circ \pi$. We sometimes write as $\pi: (X, T) \rightarrow (Y, S)$ to clarify the dynamical systems in question.

Consider a sequence of dynamical systems with factor maps connecting them:
\begin{equation*}
\xymatrix{
(X_r, T_r) \ar[r]^-{\pi_{r-1}} & (X_{r-1}, T_{r-1}) \ar[r]^-{\pi_{r-2}} & \cdots \ar[r]^-{\pi_1} & (X_1, T_1)}.
\end{equation*}
Let $d^{(i)}$ be a metric on $X_i$.

The following is Feng and Huang's original definition of weighted topological entropy \cite{Feng--Huang}. Let $\boldsymbol{w} = (w_1, w_2, \ldots, w_r)$ be a vector with $w_1 > 0$ and $w_i \geq 0$ for $i \geq 2$ (this $\boldsymbol{w}$ is the ``weight''). Let $n$ be a natural number and $\vep$ a positive number. For $x \in X_r$, define the \textbf{$\boldsymbol{n}$-th $\boldsymbol{w}$-weighted Bowen ball of radius $\boldsymbol{\vep}$ centered at $\boldsymbol{x}$} by
\begin{equation*}
B^{\boldsymbol{w}}_n(x, \vep) = \left\{ y \in X_r \spa \middle|
\begin{array}{l}
\text{$d^{(i)} \! \left( T^j_i(\pi^{(i)}(x)), T^j_i(\pi^{(i)}(y)) \right) < \vep$ for every} \\[2pt]
\text{$0 \leq j \leq \ceil{(w_1 + \cdots + w_i)n}$ and $1 \leq i \leq r$.}
\end{array}
\right\}.
\end{equation*} 
Here, 
\begin{gather*}
\pi^{(r)} = \mathrm{id}_{X_r}: X_r \to X_r, \\
\pi^{(i)} = \pi_i \circ \pi_{i+1} \circ \cdots \circ \pi_r: X_r \to X_i.
\end{gather*}
Consider $\Gamma = \{ B^{\boldsymbol{w}}_{n_j}(x_j, \vep) \}_j$, an at-most countable cover of $X_r$ by weighted Bowen balls. Let $n(\Gamma) = \min_j n_j$. For $s \geq 0$ and $N \in \mathbb{N}$, let
\[ \Lambda^{\boldsymbol{w}, s}_{N, \vep} = \inf \left\{ \sum_j e^{-sn_j} \spa \middle| \spa \text{ $\Gamma = \{ B^{\boldsymbol{w}}_{n_j}(x_j, \vep) \}_j$ covers $X_r$ and $n(\Gamma) \geq N$} \right\}. \]
This quantity is non-decreasing as $N \to \infty$. The following limit hence exists:
\[ \Lambda^{\boldsymbol{w}, s}_{\vep} = \lim_{N \to \infty} \Lambda^{\boldsymbol{w}, s}_{N, \vep}. \]
There is a value of $s$ where $\Lambda^{\boldsymbol{w}, s}_{\vep}$ jumps from $\infty$ to $0$, which we will denote by $h^{\boldsymbol{w}}_{\mathrm{wgt}}(T_r, \vep)$:
\begin{equation*}
\Lambda^{\boldsymbol{w}, s}_{\vep} = \left\{
\begin{array}{ll}
\infty & (s < h^{\boldsymbol{w}}_{\mathrm{wgt}}(T_r, \vep)) \\
0 & (s > h^{\boldsymbol{w}}_{\mathrm{wgt}}(T_r, \vep))
\end{array}
\right..
\end{equation*}
(Here, ``wgt'' stands for ``weight''.) The value $h^{\boldsymbol{w}}_{\mathrm{wgt}}(T_r, \vep)$ is non-decreasing as $\vep \to 0$. Therefore, we can define the \textbf{$\boldsymbol{w}$-weighted topological entropy} $h^{\boldsymbol{w}}_{\mathrm{wgt}}(\boldsymbol{T})$, where $\boldsymbol{T} = (T_i)_i$, by
\[ h^{\boldsymbol{w}}_{\mathrm{wgt}}(\boldsymbol{T}) = \lim_{\vep \to 0} h^{\boldsymbol{w}}_{\mathrm{wgt}}(T_r, \vep). \]

Now, we state the newly introduced definition in \cite{Alibabaei}. Take a natural number $N$ and define a new metric $d^{(i)}_N$ on $X_i$ by
\[ d^{(i)}_N(x_1, \spa x_2) = \max_{0\leq n < N} d^{(i)}({T_i}^n x_1, \spa {T_i}^n x_2). \]
Let ${\boldsymbol{a}}=(a_1, \spa a_2, \spa \cdots, a_{r-1})$ with $0 \leq a_i \leq 1$ for each $i$. Take a positive number $\vep$. We backward inductively define a quantity $\#^{\boldsymbol{a}}_i(\Omega, \spa N, \spa \vep)$ for $\Omega \subset X_i$. For $\Omega \subset X_r$, set
\begin{align*}
\#^{\boldsymbol{a}}_r(\Omega, \spa N, \spa \vep)
&= \min \left\{ \spa n \in \mathbb{N} \spa \middle|
\begin{array}{l}
\text{There exists an open cover $\{U_j\}_{j=1}^n$ of $\Omega$} \\ 
\text{with $\diam(U_j, \spa d_N^{(r)}) < \vep$ for all $ 1 \spa \leq j \spa \leq n$}
\end{array}
\right\}.
\end{align*}
(The quantity $\#^{\boldsymbol{a}}_r(\Omega, \spa N, \spa \vep)$ is independent of the parameter $\boldsymbol{a}$. However, we use this notation for the convenience of the sequel.)
Suppose $\#^{\boldsymbol{a}}_{i+1}$ is already defined. For $\Omega \subset X_{i}$, we set
\begin{flalign*}
& \#^{\boldsymbol{a}}_{i}(\Omega, \spa N, \spa \vep) &
\end{flalign*} \\[-35pt]
\begin{align*}
&= \min \left\{ \spa \sum_{j=1}^n \Big( \#^{\boldsymbol{a}}_{i+1}(\pi_i^{-1}(U_j), \spa N, \spa \vep) \Big)^{a_i} \spa \middle|
\begin{array}{l}
\text{$n \in \mathbb{N}$, $\{U_j\}_{j=1}^n$ is an open cover of $\Omega$} \\ 
\text{with $\diam(U_j, \spa d_N^{(i)}) < \vep$ for all $ 1 \spa \leq j \spa \leq n$}
\end{array}
\right\}.
\end{align*} 
We define the \textbf{topological entropy of ${\boldsymbol{a}}$-exponent} $h^{\boldsymbol{a}}_{\mathrm{exp}}(\boldsymbol{T})$, where $\boldsymbol{T} = (T_i)_i$, by \\[-1pt]
\begin{equation*}
h^{\boldsymbol{a}}_{\mathrm{exp}}(\boldsymbol{T}) = \lim_{\vep \to 0} \left( \lim_{N \to \infty} \frac{\log{\#^{\boldsymbol{a}}_1(X_1, \spa N, \spa \vep)}}{N} \right).
\end{equation*} \\[-3pt]
(Here, ``exp'' stands for ``exponent''.) This limit exists since $\log{\#^{\boldsymbol{a}}_r(X_r, \spa N, \spa \vep)}$ is sub-additive in $N$ and non-decreasing as $\vep$ tends to $0$. (We remark that the above definitions have reversed order from the original definitions. This is to simplify the following arguments.)

Via the variational principle, we obtained the following in \cite{Alibabaei}.

\begin{theorem}[{{\cite[Corollary 1.3]{Alibabaei}}}] \label{theorem: two definitions of weighted topological entropy}
For ${\boldsymbol{a}}=(a_1, \spa a_2, \spa \cdots, a_{r-1})$ with $0 < a_i \leq 1$ for each $i$,
\begin{equation*}
h^{\boldsymbol{a}}_{\mathrm{exp}}(\boldsymbol{T}) = h^{\boldsymbol{w_a}}_{\mathrm{wgt}}(\boldsymbol{T}).
\end{equation*}
Here, $\boldsymbol{w_a} = (w_1, \spa\cdots, \spa w_r)$ is defined by
\begin{eqnarray*}
\left\{
\begin{array}{l}
w_1 = a_1 a_2 a_3 \cdots a_{r-1}\\
w_2 = (1-a_1) a_2 a_3 \cdots a_{r-1} \\
w_3 = (1-a_2) a_3 \cdots a_{r-1} \\
\hspace{50pt} \vdots \\
w_{r-1} = (1-a_{r-2}) a_{r-1} \\
w_r = 1- a_{r-1}
\end{array}
\right..
\end{eqnarray*}
\end{theorem}

Theorem \ref{theorem: two definitions of weighted topological entropy} has profound implications, especially when the dynamical systems in question are affine-invariant under an endomorphism that is a direct sum of conformal endomorphisms, as in the following argument. 

\subsection{Proof of the main theorem}

Recall that $D_i = I_1 \times I_2 \times \cdots \times I_i$. The projection $p_{r-1} : S|_N \rightarrow D_{r-1}^N$ was defined by
\begin{equation*}
p_{r-1} \left( \left( e^{(n)}_1, e^{(n)}_2, \ldots, e^{(n)}_{r-1}, e^{(n)}_r \right)_{n=1}^N \right) = \left( e^{(n)}_1, e^{(n)}_2, \ldots, e^{(n)}_{r-1} \right)_{n=1}^N.
\end{equation*}
We let $S^{(r-1)} = p_{r-1}(S)$, and the rest of $p_i$ and $S^{(i)}$ were defined similarly. Let $X_i = R(S^{(i)})$ and $T_i = \mathrm{diag}(m_1, m_2, \ldots, m_i)$. Also, define factor maps $\pi_i: X_{i+1} \rightarrow X_i$ for $ 1 \leq i \leq r-1$ by $\pi_i \left( x_1, x_2, \ldots, x_i, x_{i+1} \right) = \left( x_1, x_2, \ldots, x_i \right).$ Then, we have the following chain of dynamical systems with factor maps connecting them:
\begin{equation*}
\xymatrix{
(X_r, T_r) \ar[r]^-{\pi_{r-1}} & (X_{r-1}, T_{r-1}) \ar[r]^-{\pi_{r-2}} & \cdots \ar[r]^-{\pi_{1}} & (X_1, T_1)}.
\end{equation*}

\begin{proof}[Proof of Theorem \ref{theorem: combinatorial expression of Hausdorff dimension}]

The proof is motivated by \cite[Claim 1.6]{Alibabaei}. Let
\begin{equation*}
\boldsymbol{w} = \left( \frac{\log{m_1}}{\log{m_r}}, \quad \frac{\log{m_1}}{\log{m_{r-1}}} - \frac{\log{m_1}}{\log{m_r}}, \spa \ldots \spa , \quad \frac{\log{m_1}}{\log{m_2}} - \frac{\log{m_1}}{\log{m_3}}, \quad 1 - \frac{\log{m_1}}{\log{m_2}} \right).
\end{equation*}
Then, each $n$-th $\boldsymbol{w}$-weighted Bowen ball is approximately a square of side length $\vep m_1^{-n}$. Therefore,
\begin{equation*}
\mathrm{dim}_{\mathrm{H}} (X_r) = \frac{h^{\boldsymbol{w}}_{\mathrm{wgt}}(\boldsymbol{T})}{\log{m_1}}.
\end{equation*}
Set $a_i = \log_{m_{r-i+1}} m_{r-i}$ for each $1 \leq i \leq r-1$, then $\boldsymbol{w_a}$ equals $\boldsymbol{w}$ above. We conclude from Theorem \ref{theorem: two definitions of weighted topological entropy} that
\begin{equation*}
\mathrm{dim}_{\mathrm{H}} (X_r) = \frac{h^{\boldsymbol{a}}_{\mathrm{exp}}(\boldsymbol{T})}{\log{m_1}}.
\end{equation*}
Hence, we will calculate $h^{\boldsymbol{a}}_{\mathrm{exp}}(\boldsymbol{T})$.

Fix $0 < \vep <\frac{1}{m_r}$, and take a natural number $n = n(\vep)$ with $m_1^{-n} < \vep$. Take a natural nummber $N$ and let $M = N + n$. We will later let $N$ go to infinity, and so will $M$. For $s \in S^{(i)}|_M$, define (recall that $T_i = \mathrm{diag}(m_1, m_2, \ldots, m_i)$)
\begin{equation*}
U^{(i)}_s = \left\{ \sum_{k=0}^{\infty} T_i^{-k} e_k \in X_i \spa \middle| \spa \text{$e_k \in D_i$ for each $k$ and $(e_1, \dots, e_M) = s$} \right\}.
\end{equation*}
Then, $\{U^{(i)}_s\}_{s \in S^{(i)}|_M}$ is a closed cover of $X_i$ with $\diam(U^{(i)}_x, d^{(i)}_N) < \vep$. For $s, t \in D_i^M$, we write $s \backsim t$ if and only if $U^{(i)}_s \cap U^{(i)}_t \ne \varnothing$. We have for any $1 \leq i \leq r-1$ and $s \in S^{(i)}|_M$
\[ \pi_i^{-1}(U^{(i)}_s) \subset \bigcup_{\substack{t \in S^{(i)}|_M \\ t \backsim s}} \hspace{3pt} \bigcup_{u \in {p_i}^{-1}(t)} U^{(i+1)}_u.  \]
This gives the estimation below for each $1 \leq i \leq r-2$.
\[ \#_{i+1}^{\boldsymbol{a}} \left( \pi_i^{-1}(U^{(i)}_s), \spa N, \spa \vep \right)
\leq \sum_{
		\substack{t \in S^{(i)}|_M \\ t \backsim s}}
	\sum_{u \in p_i^{-1}(t)}
	{\left( \#_{i+2}^{\boldsymbol{a}}
		\left( \pi_{i+1}^{-1}(U^{(i+1)}_u), \spa N, \spa \vep
		\right) 
	\right)}^{a_{i+1}}. \]
Also, we have
\[ \#_{r}^{\boldsymbol{a}} \left( \pi_{r-1}^{-1}(U^{(r-1)}_u), \spa N, \spa \vep \right)
\leq {\left(
		\sum_{
			\substack{v \in S^{(r-1)}|_M \\ v \backsim u}}
		\abs{p_{r-1}^{-1}(v)}
	\right)}^{a_{r-1}}. \]

Using these inequalities, we obtain
\begin{flalign*}
&
\#_{1}^{\boldsymbol{a}} \left( X_1, \spa N, \spa \vep \right) 
\begin{multlined}[t][11.5cm]
\leq
\sum_{s_1^{} \in S^{(1)}|_M} {\left( \#_{2}^{\boldsymbol{a}} \left( \pi_{1}^{-1}(U^{(1)}_{s_1^{}}), \spa N, \spa \vep \right) \right)}^{a_1}
\end{multlined} \\
&
\begin{multlined}[t][11.5cm]
\hspace{63pt}
\leq \sum_{s_1^{} \in S^{(1)}|_M}
	{\left( \sum_{
			\substack{s_1' \in S^{(1)}|_M \\ s_1' \backsim s_1^{}}} 
		\sum_{s_2^{} \in p_1^{-1}(s_1')}
		{\left( \#_{3}^{\boldsymbol{a}}
			\left( \pi_{2}^{-1}(U^{(2)}_{s_2^{}}), \spa N, \spa \vep 
			\right)
		\right)}^{a_2}
	\right)}^{a_1}
\end{multlined} &
\end{flalign*}
\begin{flalign*}
&
\begin{multlined}[t][11.5cm]
\leq \cdots \leq \sum_{s_1^{} \in S^{(1)}|_M} 
\left( \sum_{
		\substack{s_1' \in S^{(1)}|_M \\ s_1' \backsim s_1^{}}
		} \sum_{s_2^{} \in p_1^{-1}(s_1')}
\right. \\
\hspace{50pt} { \left( \cdots
			{ \left.
				{\left( \sum_{
					\substack{s_{r-2}' \in S^{(r-2)}|_M \\ s_{r-2}' \backsim s_{r-2}^{}}
					} \sum_{s_{r-1}^{} \in p_{r-2}^{-1}(s_{r-2}')}
					{\left(
						\#_{r}^{\boldsymbol{a}} 
							\left( \pi_{r-1}^{-1}(U^{(r-1)}_{s_{r-1}^{}}), \spa N, \spa \vep 
							\right)
					\right)}^{a_{r-1}}
				\right)}^{a_{r-2}} \cdots
			\right)}^{a_2}
		\right)}^{a_1}
\end{multlined} &
\end{flalign*}
\begin{flalign} \label{inequality: core inequality for Hausdorff dimension of sofic sets}
&
\begin{multlined}[t][11.5cm]
\leq \sum_{s_1^{} \in S^{(1)}|_M} 
\left( \sum_{
		\substack{s_1' \in S^{(1)}|_M \\ s_1' \backsim s_1^{}}
		} \sum_{s_2^{} \in p_1^{-1}(s_1')}
\right. \\
\hspace{50pt} { \left. \cdots 
				{\left(
					{\left( \sum_{
						\substack{s_{r-2}' \in S^{(r-2)}|_M \\ s_{r-2}' \backsim s_{r-2}^{}}} 
						\sum_{s_{r-1}' \in p_{r-2}^{-1}(s_{r-2}')}
						{\left( \sum_{
								\substack{s_{r-1}' \in S^{(r-1)}|_M \\ s_{r-1}' \backsim s_{r-1}}}
							\abs{p_{r-1}^{-1}(s_{r-1}')}
						\right)}^{a_{r-1}}
					\right)}^{a_{r-2}} \cdots
				\right)}^{a_2}
			\right)}^{a_1}.
\end{multlined} &
\end{flalign}

Let
$f^{(r-1)}(s_{r-1}) = \abs{p_{r-1}^{-1}(s_{r-1})}$ for $s_{r-1} \in S^{(r-1)}|_M$. For $1 \leq i \leq r-2$ and $s_i \in S^{(i)}|_M$, let
\begin{flalign*}
&
\begin{multlined}[t][11.5cm]
f^{(i)}(s_i^{}) = \sum_{s_{i+1}^{} \in p_i^{-1}(s_i^{})} 
\left( \sum_{
		\substack{s_{i+1}' \in S^{(i+1)}|_M \\ s_{i+1}' \backsim s_{i+1}^{}}
		} \sum_{s_{i+2}^{} \in p_{i+1}^{-1}(s_{i+1}')}
\right. \\
\hspace{50pt} { \left. \cdots 
				{\left(
					{\left( \sum_{
						\substack{s_{r-2}' \in S^{(r-2)}|_M \\ s_{r-2}' \backsim s_{r-2}^{}}} 
						\sum_{s_{r-1}' \in p_{r-2}^{-1}(s_{r-2}')}
						{\left( \sum_{
								\substack{s_{r-1}' \in S^{(r-1)}|_M \\ s_{r-1}' \backsim s_{r-1}}}
							\abs{p_{r-1}^{-1}(s_{r-1}')}
						\right)}^{a_{r-1}}
					\right)}^{a_{r-2}} \cdots
				\right)}^{a_2}
			\right)}^{a_1}.
\end{multlined} &
\end{flalign*}
Also, define for $1 \leq i \leq r-1$
\[ g^{(i)}_M =
\sum_{s_1^{} \in S^{(1)}|_M}
	{\left( \sum_{s_2^{} \in p_{1}^{-1}(s_1^{})}
		{\left( \cdots 
			{\left( \sum_{s_{i}^{} \in p_{i-1}^{-1}(s_{i-1}^{})}
				{\left( 
					\sum_{
						\substack{s_i' \in S^{(i)}|_M \\ 
								s_i' \backsim s_i^{}}}
					f^{(i)}(s_i')
				\right)}^{a_i}
			\right)}^{a_{i-1}} \cdots
		\right)}^{a_2}
	\right)}^{a_1},
\]
and
\[ g^{(r)}_M =
\sum_{s_1^{} \in S^{(1)}|_M}
	{\left( \sum_{s_2^{} \in p_{1}^{-1}(s_1^{})}
		{\left( \cdots 
			{\left( \sum_{s_{r-1}^{} \in p_{r-2}^{-1}(s_{r-2}^{})}
					\abs{p_{r-1}^{-1}(s_{r-1})}
			\right)}^{a_{r-1}} \cdots
		\right)}^{a_2}
	\right)}^{a_1}.
\]

We use induction on $i$ to show that $\#_{1}^{\boldsymbol{a}} \left( X_1, \spa N, \spa \vep \right)$ is not more than $g^{(i)}_M$ times a constant that does not depend on $M$ (thus nor on $N$). When $i=1$, it is immediate from the definition and inequality (\ref{inequality: core inequality for Hausdorff dimension of sofic sets}) that $\#_{1}^{\boldsymbol{a}} \left( X_1, \spa N, \spa \vep \right) \leq g^{(1)}_M$. For $i=2$,
\[ \#_{1}^{\boldsymbol{a}} \left( X_1, \spa N, \spa \vep \right) \leq
g^{(1)}_M =
\sum_{s_1^{} \in S^{(1)}|_M}
	{\left(\sum_{
		\substack{s_1' \in S^{(1)}|_M \\ s_1' 								\backsim s_1^{}}}
		f^{(1)}(s_1') 
	\right)}^{a_1}
\leq \sum_{s_1^{} \in S^{(1)}|_M}
	\sum_{
		\substack{s_1' \in S^{(1)}|_M \\ s_1' 								\backsim s_1^{}}}
	{\left(
		f^{(1)}(s_1') 
	\right)}^{a_1}. \]
Since the number of adjacent intervals in $\mathbb{R}$ is less than or equal to $3$, the multiplicity of $s_1'$ that appears in the sum above is finite. Thus,
\[ \#_{1}^{\boldsymbol{a}} \left( X_1, \spa N, \spa \vep \right) \leq
3 \sum_{s_1^{} \in S^{(1)}|_M}
	{\left(f^{(1)}(s_1) 
	\right)}^{a_1} = 3g^{(2)}_M. \]

Similarly, suppose that
\begin{equation} \label{equation: induction hypothesis for Hausdorff dimension of sofic sets}
\#_{1}^{\boldsymbol{a}} \left( X_1, \spa N, \spa \vep \right)
\leq 3^{\frac{i(i-1)(i+1)}{6}} g^{(i)}_M
\end{equation}
for some $i$. Notice that if $s_i' \in S^{(i)}|_M$ and $s_i' \backsim s_i^{}$, then we must have $p_{i-1}(s_i') \backsim p_{i-1}(s_i)$. Also, the number of adjacent cuboids in $\mathbb{R}^i$ is not more than $3^i$. Therefore, the following terms inside the sum of $g^{(i)}_M$ can be evaluated as
\begin{align*}
\sum_{s_{i}^{} \in p_{i-1}^{-1}(s_{i-1}^{})}
	{\left( \sum_{
			\substack{s_i' \in S^{(i)}|_M \\ s_i' \backsim s_i^{}}}
		f^{(i)}(s_i')
	\right)}^{a_i}
&\leq
\sum_{s_{i}^{} \in p_{i-1}^{-1}(s_{i-1}^{})} \hspace{3pt}
\sum_{
	\substack{s_i' \in S^{(i)}|_M \\ s_i' \backsim s_i^{}}}
	{\left(
		f^{(i)}(s_i')
	\right)}^{a_i} \\
&\leq
3^i \sum_{
	\substack{s_{i-1}' \in S^{(i-1)}|_M \\ s_{i-1}' \backsim s_{i-1}^{}}} \hspace{3pt}
\sum_{s_i' \in p_{i-1}^{-1}(s_{i-1}')}
	{\left(
		f^{(i)}(s_i')
	\right)}^{a_i}.
\end{align*}
We apply this to inequality (\ref{equation: induction hypothesis for Hausdorff dimension of sofic sets}) and continue this reasoning inductively until we have
\begin{align*}
\#_{1}^{\boldsymbol{a}} \left( X_1, \spa N, \spa \vep \right)
&\leq
3^{\frac{i(i-1)(i+1)}{6}} 3^{\frac{i(i+1)}{2}} \hspace{-2pt} \sum_{s_1^{} \in S^{(1)}|_M}
	{\left( \sum_{s_2' \in p_{1}^{-1}(s_1')}
		{ \hspace{-2pt} \left( \hspace{-2pt} \cdots 
			{\left( 
				\sum_{s_i' \in p_{i-1}^{-1}(s_{i-1}')}
				{\left(
					f^{(i)}(s_i')
				\right)}^{a_i}
			\right)}^{ \hspace{-2pt} a_{i-1}} \hspace{-4pt} \cdots
		\right)}^{ \hspace{-2pt} a_2}
	\right)}^{a_1} \\
&=
3^{\frac{i(i+1)(i+2)}{6}} g^{(i+1)}_M.
\end{align*}

Therefore, we have
\[ \#_{1}^{\boldsymbol{a}} \left( X_1, \spa N, \spa \vep \right)
\leq 3^{\frac{r(r-1)(r+1)}{6}} g^{(r)}_M. \]
This implies
\[ h^{\boldsymbol{a}}_{\mathrm{exp}}(\boldsymbol{T})
= \varliminf_{\vep \to 0} \varliminf_{N \to \infty} \frac{1}{N} \log{ \left( \#_{1}^{\boldsymbol{a}} \left( X_1, \spa N, \spa \vep \right) \right)}
\leq \varliminf_{\vep \to 0}  \varliminf_{N \to \infty} \frac{1}{N} \log{\left(g^{(r)}_{N + n(\vep)}\right)}
= \varliminf_{N \to \infty} \frac{1}{N} \log{\left(g^{(r)}_N\right)}. \]

The other direction is straightforward, as we can use $\vep$-separated sets. Fix $0 < \vep < \frac{1}{m_r}$ and a natural number $N$. We consider the following family of sets, with each element corresponding to a cuboid in ${\mathbb{T}}^i$:
\[ Q^{(i)}_N
= \left\{ C \subset D_r^{\mathbb{N}} \setcond
\begin{array}{l}
\text{There is $\left( s^{(n)} \right)_{n=1}^{N} \in S^{(i)}|_N$ such that every $\left( e^{(n)} \right)_{n=1}^{\infty} \in C$} \\[7pt]
\text{satisfies $e^{(n)} = s^{(n)}$ for each $1 \leq n \leq N$}
\end{array}
\right\}. \]
From each $C \in Q^{(i)}_N$, take an element $e_C \in C \cap S^{(i)}$ and denote by $P^{(i)}_N$ the collection of these points:
\[ P^{(i)}_N = \left\{ e_C \in C \cap S^{(i)} \setcond C \in Q^{(i)}_N \right\}. \]
Since for each $s \ne t \in P^{(i)}_N$, there is $1 \leq n \leq N$ such that the $n$-th coordinates of $s$ and $t$ are different (otherwise, $s =t$), we see that $Y^{(i)}_N = R(P^{(i)}_N)$ is an $\vep$-separated set with respect to the distance $d^{(i)}_N$ on $X_i = R(S^{(i)})$.

Next, consider an arbitrary chain of open $(N, \vep)$-covers of $(X_i, T_i)_{i=1}^r$, defined as follows.
\begin{enumerate}
\item[(1)] For every $i$ and $V \in \flo{F}^{(i)}$, we have $\diam(V, d^{(i)}_N) < \vep$.
\item[(2)] For each $1 \leq i \leq r-1$ and $U \in \flo{F}^{(i)}$, there is $\flo{F}^{(i+1)}(U) \subset \flo{F}^{(i+1)}$ such that
\[ \pi_{i}^{-1}(U) \subset \bigcup \flo{F}^{(i+1)}(U) \]
and
\[ \flo{F}^{(i+1)} = \bigcup_{U \in \flo{F}^{(i)}} \flo{F}^{(i+1)}(U). \]
\end{enumerate}
(The above definition appears in \cite[Definition 3.1]{Alibabaei}. Also, \cite[Remark 3.2]{Alibabaei} is relevant.)

Now, by condition (1), we have $\abs{V \cap Y^{(i)}_N} \leq 1$ for each $V \in \flo{F}^{(i)}$. Suppose that for some $V \in \flo{F}^{(i)}$, we have $R(s_{i}^{}) \in V \cap Y^{(i)}_N$ with $s_{i}^{} = \left( s_1^{(n)}, \ldots, s_{i}^{(n)} \right) \in P^{(i)}_N$. Then $\pi_{i}^{-1}(V)$ must contain at least $\abs{p_{i}^{-1}(s_{i}^{})}$ elements of $Y^{(i+1)}_N$. Since these points are $\vep$-separated, this in turn implies that $\flo{F}^{(i+1)}(V)$ from condition (2) has at least $\abs{p_{i}^{-1}(s_{i}^{})}$ elements. Thus,

\[ \#_{1}^{\boldsymbol{a}} \left( X_1, \spa N, \spa \vep \right)
\geq \sum_{s_1^{} \in S^{(1)}|_M}
	{\left( \sum_{s_2^{} \in p_{1}^{-1}(s_1^{})}
		{\left( \cdots 
			{\left( \sum_{s_{r-1}^{} \in p_{r-2}^{-1}(s_{r-2}^{})}
					\abs{p_{r-1}^{-1}(s_{r-1})}
			\right)}^{a_{r-1}} \cdots
		\right)}^{a_2}
	\right)}^{a_1}
= g^{(r)}_N. \]
This yields the desired inequality;
\[ h^{\boldsymbol{a}}_{\mathrm{exp}}(\boldsymbol{T})
= \varlimsup_{\vep \to 0} \varlimsup_{N \to \infty} \frac{1}{N} \log{ \left( \#_{1}^{\boldsymbol{a}} \left( X_1, \spa N, \spa \vep \right) \right)}
\geq \varlimsup_{\vep \to 0}  \varlimsup_{N \to \infty} \frac{1}{N} \log{\left(g^{(r)}_{N}\right)}
= \varlimsup_{N \to \infty} \frac{1}{N} \log{\left(g^{(r)}_N\right)}. \]

We conclude that
\[ h^{\boldsymbol{a}}_{\mathrm{exp}}(\boldsymbol{T})
= \lim_{N \to \infty} \spa \frac{1}{N} \spa \log{\hspace{-3pt}\sum_{s_1^{} \in S^{(1)}|_N} \hspace{-3pt}
\left( \sum_{s_2^{} \in p_1^{-1}(s_1^{})} \hspace{-5pt}
	\left( \cdots
		\left( \sum_{s_{r-1}^{} \in p_{r-2}^{-1}(s_{r-2}^{})} \hspace{-12pt}
		{\abs{p_{r-1}^{-1}(s_{r-1}^{})}}^{a_{r-1}}
		\right)^{a_{r-2}} \hspace{-9pt} \cdots
	\right)^{a_2}
\right)^{a_1}}. \]
\end{proof}

\vspace{0.5cm}

\address{
Department of Mathematics, Kyoto University, Kyoto 606-8501, Japan}

\textit{E-mail}: \texttt{alibabaei.nima.28c@st.kyoto-u.ac.jp}

\end{document}